\documentclass[12pt]{amsart}
\usepackage{amssymb}
\usepackage{times}
\usepackage{graphicx}

\theoremstyle{plain}

\numberwithin{equation}{section}

\newcommand{\calB}{\mathcal{B}}

\newcommand{\calD}{\mathcal{D}}

\newcommand{\calH}{\mathcal{H}}

\newcommand{\calM}{\mathcal{M}}

\newcommand{\bbF}{\mathbb{F}}
\newcommand{\bbC}{\mathbb{C}}

\newcommand{\bbP}{\mathbb{P}}
\newcommand{\bbQ}{\mathbb{Q}}
\newcommand{\bbR}{\mathbb{R}}

\newcommand{\bbZ}{\mathbb{Z}}

\def\O{{\text{O}}}

\def\PGL{{\text{PGL}}}

\def\rank{{\text{rank}}}
\def\det{{\text{det}}}
\def\mod{{\text{mod}}}

\begin{document}
\title [The moduli space of 5 points on $\bbP^1$] {The moduli space of 5 points on $\bbP^1$ and K3 surfaces}
\author{Shigeyuki Kond{$\bar{\rm o}$}}
\address{Graduate School of Mathematics, Nagoya University, Nagoya,
464-8602, Japan}
\email{kondo@math.nagoya-u.ac.jp}
\thanks{Research of the author is partially supported by
Grant-in-Aid for Scientific Research A-14204001, Japan}

\begin{abstract}
We show that the moduli space of ordered 5 points on $\bbP^1$ is isomorphic to an arithmetic quotient
of a complex ball by using the theory of periods of $K3$ surfaces.  We also discuss a relation between our
uniformization and the one given by Shimura \cite{S}, Terada \cite{Te}, Deligne-Mostow \cite{DM}.
\end{abstract}
\maketitle

\rightline{\it Dedicated to Professor Yukihiko Namikawa on his 60th birthday}

\section{Introduction}
The purpose of this note is to show that the moduli space of 5 ordered points on $\bbP^1$
is isomorphic to an arithmetic quotient of a 2-dimensional complex ball by using the theory of
periods of $K3$ surfaces (Theorem \ref{Main}).  This was announced in \cite{K2}, Remark 6.  
The main idea is to associate a $K3$ surface with an automorphism
of order 5 to a set of 5 ordered points on  $\bbP^1$ (see  \S \ref{assK3}).  The period domain of such $K3$ surfaces
is a 10-dimensional bounded symmetric domain of type $IV$.  
We remark  that a non-zero holomorphic 2-form on the $K3$ surface is an
eigen-vector of the automorphism, which implies that 
the period domain of the pairs of these $K3$ surfaces and the automorphism 
of order 5 is a 2-dimensional complex ball associated to a hermitian form of the signature $(1,2)$ defined over $\bbZ[\zeta]$
where $\zeta$ is a primitive 5-th root of unity (see \ref{period}).  Here we use several fundamental results of Nikulin
\cite{N1}, \cite{N2}, \cite{N3} on automorphisms of $K3$ surfaces and the lattice theory.

Note that this moduli space is isomorphic to
the moduli space of nodal del Pezzo surfaces of degree 4.  For  the moduli space of 
del Pezzo surfaces of degree 1, 2 or  3, the similar description holds. 
See \cite{K2}, Remark 5, \cite{K1}, \cite{DGK}, respectively.  

On the other hand, Shimura \cite{S}, Terada \cite{Te}, Deligne-Mostow \cite{DM} gave a complex ball
uniformization by using the periods of the curve $C$ which is the 5-fold cyclic covering of $\bbP^1$ branched along 5 points.  
We shall discuss a relation between their uniformization and ours in  \S  \ref{Shimura}.
In fact, the above $K3$ surface has an isotrivial pencil whose general member is the unique 
smooth curve $D$ of genus 2 admitting an automorphism of order 5
(see Lemma \ref{pencil}).  We show that the above $K3$ surface is birational to the quotient of $C \times D$ by a diagonal
action of $\bbZ/5\bbZ$ in \S \ref{Shimura}.

In this paper, a {\it lattice} means a $\bbZ$-valued non-degenerate symmetric bilinear form on 
a free $\bbZ$-module of finite rank.  We denote by $U$ or $V$ the even lattice defined by the matrix
$
\begin{pmatrix}0&1
\\1&0
\end{pmatrix}
$,
$
\begin{pmatrix}2&1
\\1&-2
\end{pmatrix}
$,
respectively and by $A_m$, $D_n$ or $E_l$ the even negative definite lattice defined by the Dynkin matrix of
type $A_m$, $D_n$ or $E_l$ respectively.  If $L$ is a lattice and $m$ is an integer, 
we denote by $L(m)$ the lattice over the same $\bbZ$-module with the symmetric bilinear form
multiplied by $m$.  We also denote by $L^{\oplus m}$ the orthogonal direct sum of $m$ copies
of $L$, by $L^*$ the dual of $L$ and by $A_L$ the finite abelian group $L^*/L$.

\medskip
\noindent
{\bf Acknowledgments.} The author thanks to Igor Dolgachev for stimulating discussions and useful suggestions.  In particular the result in \ref{5-fold} is due to him. 

\section{Quartic Del Pezzo surfaces}\label{}

\subsection{Five points on $\bbP^1$}\label{}
Consider the diagonal action of $\PGL(2)$ on $(\bbP^1)^5$.  In this case, the semi-stable
points and stable points in the sense of \cite{Mu} 
coincide and the geometric quotient $P_1^5$ is smooth and compact.
The stable points are $\{p_1,...,p_5\}$ no three of which coincide.
It is known that $P_1^5$ is isomorphic to the quintic del Pezzo surface $\calD_5$, that is, a
smooth surface obtained by blowing up four points 
$\{q_1,...,q_4\}$ in general position on $\bbP^2$ (e.g. Dolgachev \cite{D}, Example 11.5).
The quintic del Pezzo surface $\calD_5$ contains 10 lines corresponding to the
4 exceptional curves over $q_1,...,q_4$ and the proper transforms of 6 lines through two
points from $\{q_1,...,q_4\}$.  These ten lines correspond to the locus consisting of
$\{p_1,...,p_5\}$ with $p_i = p_j$ for some $i, j$.
The group of automorphisms of $\calD_5$ is
isomorphic to the Weyl group $W(A_4) \simeq S_5$ which is induced from the natural action of $S_5$ on
$(\bbP^1)^5$.

\subsection{Quartic Del Pezzo surfaces}\label{}
Let $S$ be a smooth quartic Del Pezzo surface.  It is known that
$S$ is a complete intersection of two quadrics in $\bbP^4$.
Consider the pencil of quadrics whose base locus is $S$.
Its discriminant is a union of distinct five points of $\bbP^1$.  Conversely
any distinct five points $(1:\lambda_i)$ on $\bbP^1$,
the intersection of quadrics
\begin{equation}\label{quartic}
\sum_{i=1}^5 z_i^2 = \sum_{i=1}^5 \lambda_i z_i^2 = 0
\end{equation}
is a smooth quartic Del Pezzo surface.
Thus the moduli space of smooth quartic Del Pezzo surfaces is
isomprphic to $(((\bbP^1)^5 \setminus \Delta)/\PGL(2))/S_5$ where
$\Delta$ is the locus consisting of points $(x_1,...,x_5)$ with
$x_i = x_j$ for some $i, j$ \ ($i \not= j$).  
If five points are not distinct, but stable, the equation \eqref{quartic} defines
a quartic Del Pezzo surface with a node.  Thus $P_1^5$ is the coarse moduli space of
nodal quartic Del Pezzo surfaces.

\section{K3 surfaces associated to five points on $\bbP^1$}\label{assK3}

\subsection{A plane quintic curve}\label{quintic}

Let $\{ p_1,...,p_5\}$ be an ordered stable point in $(\bbP^1)^5$.
It defines a homogenious polynomial $f_5(x_1,x_2)$ of degree 5.
Let $C$ be the plane quintic curve defined by
\begin{equation}\label{quintic}
x_{0}^{5} = f_5(x_1,x_2) = \prod_{i=1}^5 (x_1 - \lambda_i x_2)
\end{equation}
The projective transformation
\begin{equation}\label{auto}
g : (x_{0} : x_{1} : x_{2})
\longrightarrow (\zeta x_{0} : x_{1} : x_{2})
\end{equation}
acts on $C$ as an automorphism of $C$ of order 5 where $\zeta$
is a primitive 5-th root of unity.
Let $E_{0}$, $L_{i}$ ($1 \leq i \leq 5$)
be lines defined by
$$E_{0} : x_{0} = 0,$$
$$L_{i}: x_{1} = \lambda_{i} x_{2}.$$
Note that all $L_{i}$ are members of the pencil of lines through
$(1:0:0)$
and $L_{i}$ meets $C$ at $(0:\lambda_{i}:1)$ with multiplicity 5.

\subsection{K3 surfaces}\label{K3}
Let $X$ be the minimal resolution of the double
cover of $\bbP^{2}$ branched
along the sextic curve $E_{0} + C$.  Then $X$ is a $K3$ surface.  
We denote by $\tau$ the covering transformation.  The projective
transformation $g$ in \eqref{auto} induces an automorphism $\sigma$ of $X$ of order 5.
We denote by the same symbol $E_{0}$
the inverse image of $E_{0}$.

\medskip
{\it Case} (i)  
Assume that the equation $f_5 = 0$ has no multiple roots.  In this case there are
5 $(-2)$-curves, denoted by $E_{i}$ ($1 \leq i \leq 5$),
obtained as exceptional curves of the minimal resolution
of singularities of type $A_1$ corresponding to the intersection of $C$ and $E_0$. 
The inverse image of $L_{i}$
is the union of two smooth rational curves $F_{i}, G_{i}$ such that
$F_{i}$ is tangent to $G_{i}$ at one point. 
Let $p, q$ be the inverse image of $(1:0:0)$.
We may assume that all $F_{i}$ ( resp. $G_{i}$ )
are through $p$ (resp. $q$).  Obviously $\sigma$ preserves each curve $E_i, F_j, G_j$ $(0\leq i \leq 5, 1\leq j \leq 5)$
and $\tau$ preserves each $E_i$ and $\tau (F_i) = G_i$.

\medskip
{\it Case} (ii)  
If $f_5 = 0$ has a multiple root, then the double cover has a rational double point 
of type $D_7$.  Hence $X$ contains 7 smooth rational curves $E_j'$, $(1 \leq j \leq 7)$
whose dual graph is of type $D_7$.  We assume that $E_1'$ meets $E_0$ and
$\langle E_1', E_2' \rangle = \langle E_2', E_3' \rangle = 
\langle E_3', E_4' \rangle = \langle E_4', E_5' \rangle =
\langle E_5', E_6' \rangle = \langle E_5', E_7' \rangle = 1$.
If $\lambda_{i}$ is a multiple root, then $F_i$ and $G_i$ are
disjoint  and each of them meets one componet of $D_7$, for example, 
$F_i$ meets $E_6'$ and $G_i$ meets $E_7'$.

\subsection{A pencil of curves of genus two }\label{}
The pencil of lines on $\bbP^{2}$ through $(1:0:0)$
gives a pencil of curves of genus two on $X$.  Each member of this pencil is
invariant under the action of the automorphism $\sigma$ of order 5.  Hence
a general member is a  smooth curve of genus two with an
automorphism of order five.  Such a curve is unique up to isomorphism
and is given by
\begin{equation}\label{bolza}
y^{2} = x(x^{5} + 1)
\end{equation}
\noindent
 (see Bolza \cite{Bol}).  If  $\lambda_{i}$ is a simple root of
the equation $f_5 = 0$, then the line $L_i$ defines a singular member of
this pencil consisting of three smooth rational curves 
$E_{i} + F_{i} + G_{i}$. 
We call this singular member a {\it singular member of type} I.
If $\lambda_{i}$ is a multiple root of $f_5 = 0$, then 
the line $L_i$ defines a singular member 
consisting of nine smooth rational curves $E_1',..., E_7', F_i, G_i$.
We call this a {\it singular member of type} II.  
The two points $p, q$ are the base points of the pencil.  After blowing up at $p, q$, 
we have a base point free pencil of curves of genus two.
The singular fibers of such pencils are completly classified by Namikawa and Ueno
\cite{NU}.  The type I (resp. type II)
corresponds to [IX-2] (resp. [IX-4] ) in  \cite{NU}.  We now conclude:

\subsection{Lemma}\label{pencil}
{\it
The pencil of lines on $\bbP^{2}$ through $(1:0:0)$
gives a pencil of curves of genus two on $X$.
A general member is a  smooth curve of genus two with an
automorphism of order five.  In case that $f_5 = 0$ has no
multiple roots, it has five singular members of type $I$.
In case that $f_5 = 0$ has a multiple root $($resp. two multiple roots$)$,
it has three singular members of type $I$ and one singular member 
of type $II$ $($resp. one of type $I$ and two of type $II$$)$.
}

\subsection{ A 5-fold cyclic cover of $\bbP^1 \times \bbP^1$}\label{5-fold}
The following is due to I. Dolgachev.
The above $K3$ surface has an automorphism 
of order 5 by construction.  This implies that $X$ is obtained as a 5-fold cyclic cover of
a rational surface.  Let $D$ be a divisor of $\bbP^1 \times \bbP^1$ defined by

\begin{equation}\label{dolgachev}
D = 4(l_1 + \cdot \cdot \cdot + l_5) + m_1 + m_2 + 3m_3
\end{equation}
where $l_1,..., l_5$ are the fibers of the first projection from $\bbP^1 \times \bbP^1$
over the five points determined by the polynomial $f_5(x_1,x_2)$ in \eqref{quintic}, and
$m_1, m_2, m_3$ are three fibers of the second projection which are unique up to
projective transformations.  Take the 5-cyclic cover of $\bbP^1 \times \bbP^1$ branched
along $D$.  Then taking the normalization and resolving the singularities, and blowing down
the proper transforms of $l_1,..., l_5, m_1, m_2, m_3$ which are exceptional curves of the first kind,
we have a $K3$ surface $Y$.   Locally the singularities over the intersection points of $l_i$ and
$m_1, m_2$ are  given by $z^5 = x^4y$ and those over the intersection points of $l_i$ and $m_3$ are given by $z^5 = x^4y^3$.
Note that the ruling of the first projection from
$\bbP^1 \times \bbP^1$ gives a pencil of curves of genus 2 on $Y$.
On the other hand, consider the involution of $\bbP^1 \times \bbP^1$ which changes
$m_1$ and $m_2$, and fixes $m_3$.  Let $m_4$ be the another fixed fiber of this
involution.  Then this involution induces an involution of $Y$ which fixes
the inverse image $C$ of $m_4$.  The last curve $C$ corresponds to the plane quintic curve
given in \eqref{quintic}.

\section{Picard and transcendental lattices}\label{Pic}

In this section we shall study the Picard lattice and the transcendental lattice of
$K3$ surfaces $X$ given in \ref{K3}.  We denote by $S_X$ the Picard lattice of $X$ and
by $T_X$ the transcendental lattice of $X$.

\subsection{The Picard lattice}\label{}

\subsection{Lemma}\label{lattice}
{\it
Assume that $f_5 = 0$ has no multiple roots.
Let $S$ be the sublattice generated by $E_0$ and components of
the singular members of the pencil in Lemma \ref{pencil}.
Then  $\rank (S) = 10$ and  $\det(S) = 5^{3}$.  Moreover if $X$ is generic
in the sense of moduli, then the Picard lattice $S_X = S$.
}

\begin{proof}
First note that the dimension of $P_1^5$ is 2.  On the other hand, 
$X$ has an automorphism $\sigma$ of order 5 induced from $g$ given in
\eqref{auto} which acts non trivially on $H^0(X, \Omega^2)$.  Nowhere vanishing
holomorphic 2-forms are eigenvectors of $\sigma^*$.
We can see that the dimension of the period domain is $(22 - \rank(S_X))/(5-1)$
(\cite{N2}, Theorem 3.1.  Also see the following section \ref{uniformization}).
Hence the local Torelli theorem implies that $\rank(S_X) = 10$ for generic $X$.
Let $S_{0}$ be the sublattice of $S_X$ 
generated by $E_{i}, F_{i}, (1 \leq i \leq 5)$.  Then a direct
calculation shows that
$\rank (S_{0}) = 10$ and $\det(S_{0}) = \pm 5^{5}$.
The first assertion now follows from the relations:
$$5E_{0} = \sum_{i=1}^{5} (F_{i} - 2 E_{i}),$$
$$G_{i} + F_{i} = 2E_{0} + \sum_{j\not= 0,i} E_{j}.$$
Note that $S^*/S \simeq ({\bbZ}/5{\bbZ})^3$.  
Now assume that $\rank(S_X) = 10$.  If $S_X \not= S$, then $S \subset S_X \subset S^*$ 
and hence there exists an algebraic cycle $C$ not contained in $S$ and satisfying
$$5C = \sum_{i=0}^{5} a_i E_i + \sum_{i=1}^{4} b_i F_i,  \ a_i, b_i  \in  \bbZ.$$
By using the relations
$$\langle 5C, E_i \rangle \equiv 0 \ (\mod \ 5), \quad 
\langle 5C, F_i \rangle \equiv 0 \ (\mod \ 5),$$
we can easily show that 
$$a_i \equiv 0 \ (\mod \ 5), \quad b_i \equiv 0 \ (\mod \ 5).$$
This is a contradiction.  
\end{proof}

\subsection{Discriminant quadratic forms}\label{}
Let $L$ be an even lattice.  We denote by $L^*$ the dual of $L$ and put $A_L = L^*/L$.
Let $$q_{L}  : A_L  \to {\bbQ}/2{\bbZ}$$ be the discriminant quadratic form defined by
$$q_L(x \ \rm{ mod} \ L) = \langle x, x \rangle  \ \rm{ mod} \ 2{\bbZ}$$
and 
$$b_L : A_L \times A_L \to {\bbQ}/{\bbZ}$$
the discriminant bilinear form defined by
$$b_L(x \ \rm{ mod} \ L, y \ \rm{ mod} \ L) = \langle x, y \rangle  \ \rm{ mod} \ {\bbZ}.$$
Let $S$ be as in Lemma \ref{lattice}.  Then
$A_{S}$ is generated by 
$$\alpha = (E_1 + 2F_1 + 3F_2 + 4E_2)/5,\  \beta = (E_1 + 2F_1 + 3F_3 + 4E_3)/5,\ \gamma = (E_1 + 2F_1 + 3F_4 + 4E_4)/5$$
with $q_{S}(\alpha) = q_{S}(\beta) = q_{S}(\gamma) = -4/5$, and 
$b_{S}(\alpha, \beta) =b_{S}(\beta, \gamma) =b_{S}(\gamma, \alpha) = 3/5$.

\subsection{The transcendental lattice}\label{}

Let $T$ be the orthogonal complement of $S$ in $H^2(X, {\bbZ})$.  For generic $X$, $T$ is isomorphic to the
transcendental lattice $T_X$ of $X$ which
consists of transcendental cycles, that is, cycles not perpendicular to holomorphic 2-forms on $X$.

\subsection{Lemma}\label{lattices}
{\it
Assume that $f_5$ has no multiple roots.  Then
$$S \simeq V \oplus A_4 \oplus A_4, \
\quad
T \simeq U \oplus V \oplus A_4 \oplus A_4
$$
where $V$ or $U$ is the lattice defined by the matrix
$\begin{pmatrix}2&1\\1&-2\end{pmatrix}$, \
$\begin{pmatrix}0&1\\1&0\end{pmatrix}
$,
respectively.
}

\begin{proof}
We can see that  $q_S$ and the discriminant quadratic form of $V \oplus A_4 \oplus A_4$ coincide.  
Also note that $q_T = -q_S$ \ (Nikulin \cite{N1}, Corollary 1.6.2).  
Now the assertion follows from Nikulin \cite{N1}, Theorem 1.14.2.
\end{proof}

\subsection{Lemma}\label{degenerate}
{\it
Let $S_i$ be the sublattice generated by $E_0$ and components of
the singular members of the pencil in Lemma \ref{pencil} where $i = 1$ or $2$ is the
number of multiple roots of $f_5 = 0$.  Let $T_i$ be 
the orthogonal complement of $S_i$ in $H^{2}(X, {\bbZ})$.  Then }
$$S_1 \simeq V \oplus E_8 \oplus A_4,  \ 
\quad
T_1 \simeq U \oplus V \oplus  A_4,
$$
$$S_2 \simeq V \oplus E_8 \oplus E_8,  \ 
\quad
T_2 \simeq U \oplus V.
$$

\begin{proof}
The proof is similar to those of Lemmas \ref{lattice}, \ref{lattices}
\end{proof}

\subsection{The K\"ahler cone}\label{Kahler}

Let $S_X$ be the Picard lattice of $X$.  Denote by $P(X)^+$ the connected component of
the set $\{ x \in S_X\otimes {\bbR} : \langle x, x \rangle > 0 \}$ which contains an ample class.
Let $\Delta(X)$ be the set of effective classes $r$ with $r^2 = -2$.
Let
$$C(X) = \{ x \in P(X)^+ : \langle x, r \rangle > 0, r \in \Delta(X) \}$$
which is called the K\"ahler cone of $X$.
It is known that $C(X) \cap S_X$ consists of ample classes.
Let $W(X)$ be the subgroup of $O(S_X)$ generated by reflections defined by
$$s_r : x \to x + \langle x, r \rangle r, \quad r\in \Delta(X).$$
Note that the action of $W(X)$ on $S_X$ can be extended to $H^2(X, {\bbZ})$ acting trivially on $T_X$
because $r \in S_X = T_X^{\perp}$.  The K\"ahler cone $C(X)$ is a fundamental domain of the action of
$W(X)$ on $P(X)^+$.

\section{Automorphisms}\label{Auto}

We use the same notation as in \S 3, 4.
In this section we study the covering involution $\tau$ of $X$ over $\bbP^2$ and the automorphism
$\sigma$ of $X$ of order 5.   

\subsection{The automorphism of order 2}\label{}

\subsection{Lemma}\label{inv}
{\it
Let $\iota = \tau^*$.  Then
the invariant sublattice $M = H^{2}(X, {\bbZ})^{<\iota>}$ is generated by
$E_{i}$ $(0 \leq i \leq 5).$
}

\begin{proof}
Note that $M$ is a 2-elementary lattice, that is, its discriminant group $A_M = M^*/M$ is a finite 2-elementary abelian group.
Let $r$ be the rank of  $M$ and let
$l$ be the number of minimal generator of $A_M \cong (\bbZ/2\bbZ)^l$.
The set of fixed points of $\tau$ is the union of $C$ and $E_{0}$.
It follows from Nikulin \cite{N3}, Theorem 4.2.2 that
$(22 - r - l)/2 = g(C) = 6$ and the number of components of fixed points set of $\tau$ other than $C$ is 
$(r - l)/2 = 1$.  Hence $r = 6, l = 4$.
On the other hand we can easily see that $\{ E_{i} : 0 \leq i \leq 5 \}$
generates a sublattice of  $M$ with rank 6 and
discriminant $2^{4}$.  Now the assertion follows.
\end{proof}

\subsection{Lemma}\label{}
{\it
Let $N$ be the orthogonal complement of 
$M = H^{2}(X, \bbZ)^{<\iota>}$ in $S$.  Then $N$ is generated by the classes of
$F_i - G_i$ $(1 \leq i \leq 5)$ and contains no $(-2)$-vectors.
}

\begin{proof}
Since $\tau (F_i) = G_i$, the classes of $F_i - G_i$ are contained in $N$.  
A direct calculation shows that their intersection matrix 
$(\langle F_i-G_i, F_j-G_j\rangle )_{1\leq i,j \leq 4}$ is

$$\begin{pmatrix}-8&2&2&2\\2&-8&2&2\\2&2&-8&2\\2&2&2&-8\end{pmatrix}$$
\noindent
whose discriminant is $\pm 2^4\cdot 5^3$.
On the other hand, $N$ is the orthogonal complement of $M$ in $S$, and $M$ (resp . $S$) has the discriminant
$\pm 2^4$ (resp. $\pm 5^3$).  Hence the discriminant of $N$ is $\pm 2^4\cdot 5^3$.
Therefore the first assertion follows.   It follows from the above intersection matrix that $N$ contains no
$(-2)$-vectors. 
\end{proof}

\subsection{Lemma}\label{-2}
{\it
Let $r$ be a $(-2)$-vector in $H^{2}(X, {\bbZ})$.  Assume that $r \in M^{\perp}$ in $H^2(X, {\bbZ})$.
Then $\langle r, \omega_X \rangle \not= 0$.
}

\begin{proof}
Assume that $\langle r, \omega_X \rangle = 0$.  Then $r$ is represented by a divisor.
By Riemann-Roch theorem, we may assume that $r$ is effective.  By assumption $\iota(r) = -r$.
On the other hand the automorphism preserves effective divisors, which is a contradiction.
\end{proof}

\subsection{Lemma}\label{V}
{\it
Let $P(M)^+$ be the connected component of the set 
$$\{ x \in M\otimes \bbR : \langle x, x \rangle > 0 \}$$
which contains the class of $C$ where $C$ is the fixed curve of $\tau$ of genus $6$.
Put 
$$C(M) = \{ x \in P(M)^+ : \langle x, E_i \rangle > 0, i = 0,1,...,5 \}.$$ 
Let $W(M)$ be the subgroup generated by
reflections associated with $(-2)$-vectors in $M$.  Then $C(M)$ is a fundamental domain of the action of
$W(M)$ on $P(M)^+$ and $O(M)/\{\pm 1\}\cdot W(M) \cong S_5$ where $S_5$ is the 
symmetry group of degree $5$ which is the automorphism group of $C(M)$.
}

\begin{proof}
First consider the dual graph of $E_0,...,E_5$.  Note that any maximal extended Dynkin diagram in this
dual graph is $\tilde{D}_4$ which has the maximal rank 4 ($ = {\rm rank}(M)-2$).
It follows from Vinberg \cite{V}, Theorem 2.6 that the group $W(M)$
is of finite index in $O(M)$ the orthogonal group of $M$.  The assertion now follows from Vinberg \cite{V}, Lemma 2.4.
\end{proof}

\subsection{Lemma}\label{V2}
{\it
Let $\tilde{W}(M)$ be the subgroup of $O(M)$ generated by
all reflections associated with negative norm vectors in $M$.  
Then $C(M)$ is a fundamental domain of the action of
$\tilde{W}(M)$ on $P(M)^+$.  Moreover $\tilde{W}(M) = W(M)\cdot S_5$.
}

\begin{proof}
First note that $W(M) \subset \tilde{W}(M) \subset O(M)$.
Let $r = E_i - E_j$ $(1\leq i < j \leq 5)$ which is a $(-4)$-vector in $M$.   Since $\langle E_i - E_j, M\rangle \subset 2\bbZ$,
the reflection defined by
$$s_r:  x \to x + \langle x, r \rangle r/2$$ 
is contained in $O(M)$.  These reflections
generate $S_5$ acting on $C(M)$ as the automorphism group of $C(M)$.  Now
Lemma \ref{V} implies that $O(M) = \{ \pm 1\}\cdot \tilde{W}(M)$.
\end{proof}

\subsection{Lemma}\label{V3}
{\it
Let $C(X)$ be the K\"ahler cone of $X$.  Then
$$C(M) = C(X) \cap P(M)^+.$$
}

\begin{proof}
Since the class of $C$ is contained in the closure of $C(X)$, 
$C(X) \cap P(M)^+ \subset C(M)$, and hence it suffices to see that any face of $C(X)$ does not cut
$C(M)$ along proper interior points of $C(M)$.
Let $r$ be the class of an effective cycle with $r^2 = -2$.  If $r \in M$, Lemma  \ref{V} implies the assertion.
Now assume $\iota(r) \not= r$.  Then $r = (r + \iota(r))/2 + (r- \iota(r))/2$.  By Hodge index theorem,
$(r- \iota(r))^2 < 0$.  Since $r^2 = -2$, this implies that $((r + \iota(r))/2)^2 \geq 0$ or $-1$.
If $((r + \iota(r))/2)^2 \geq 0$, again by Hodge index theorem, $\langle x, r + \iota(r) \rangle > 0$ for any $x \in C(M)$.
Since $\iota$ acts trivially on $M$, we have $\langle x, r \rangle > 0$.
If $((r + \iota(r))/2)^2 = -1$, then $\langle r, \iota(r) \rangle = 0$.  
Note that for any $x \in M$, $\langle x, r+\iota(r) \rangle = 2\langle x, r\rangle \in 2\bbZ$.  Hence
the $(-4)$-vector $r + \iota(r)$ defines a reflection in $\tilde{W}(M)$.
It follows from
Lemma \ref{V2} that $\langle r + \iota(r), x \rangle > 0$ for any $x \in C(M)$.
Since $\iota$ acts trivially on $M$, $\langle r, x \rangle > 0$ for any $x \in C(M)$.
Thus we have proved the assertion.
\end{proof}

\subsection{An isometry of order five}\label{isometry}
Let $\sigma$ be the automorphism of $X$ of order 5 induced by the automorphism given in \ref{auto}.
In the following Lemma \ref{isometry5}
we shall show that $\sigma^{*} \mid T$ is conjugate to the isometry $\rho$
defined as follows:

Let $e, f$ be a basis of 
$U =\begin{pmatrix}0&1\\1&0\end{pmatrix}$
satisfying $e^2 = f^2 = 0, \langle e, f \rangle = 1$.
Let $x, y$ be a basis of 
$V = \begin{pmatrix}2&1\\1&-2\end{pmatrix}$ 
satisfying $x^2 = -y^2 = 2,  \langle x, y \rangle = 1$,
and let $e_1, e_2, e_3, e_4$ be a basis of $A_4$
so that $e_i^2 = -2, \langle e_i, e_{i+1} \rangle = 1$
and other $e_i$ and $e_j$ are orthogonal.

Let $\rho_0$ be an isometry of 
$U \oplus V$ defined by 

\begin{equation}\label{}
\rho_0(e) = -f, \quad
\rho_0(f) = -e-f-y, \quad
\rho_0(x) = f-x, \quad
\rho_0(y) = 3f-x+y.
\end{equation}
\noindent
Also let $\rho_4$ be an isometry of $A_4$ defined by
\begin{equation}\label{A4}
\rho_4(e_1)=e_2, \quad
\rho_4(e_2) = e_3, \quad
\rho_4(e_3) = e_4, \quad
\rho_4(e_4) 
= -(e_1+e_2+e_3+e_4).
\end{equation}
\noindent
Combaining $\rho_0$ and $\rho_4$, we define an isometry $\rho$ of $T= U\oplus V \oplus A_4\oplus A_4$.
By definition, $\rho$ is of order 5 and has no non-zero fixed vectors in $T$.
Moreover the action of $\rho$ on the discriminant group $T^*/T$ is trivial.
Hence $\rho$ can be extended to an isometry $\rho$ (we use the same
symbol) of $H^{2}(X, {\bbZ})$ acting trivially on $S$ (Nikulin \cite{N1}, Corollary 1.5.2).

\subsection{Lemma}\label{isometry5}
{\it
The isometry $\sigma^*$ is conjugate to $\rho$.
}

\begin{proof}
By the surjectivity of the period map of $K3$ surfaces,
there exists a $K3$ surface $X'$ whose transcendental lattice $T_{X'}$ is isomorphic to $T$.  Moreover
we may assume that $\omega_{X'}$ is an eigenvector of $\rho$ under the isomorphism $T_{X'} \cong T$. 
Since $\rho$ acts trivially on $S_{X'}$, there exists an automorphism $\sigma'$ of $X'$ with $(\sigma')^* = \rho$
(\cite{PS}).  

Since $S_{X'} \cong S$, there exist 16 $(-2)$-classes in $S_{X'}$ whose dual graph coincides with that of
$E_i$, $(0 \leq i \leq 5)$,  $F_j, G_k$ $(1\leq j, k \leq 5)$ on $X$ in \ref{K3}.  
We denote by $E_i', F_j', G_k'$ these divisors corresponding to $E_i, F_j, G_k$.  We shall show that 
if necessary by changing them by $w(E_i'), w(F_j'), w(G_k')$ for a suitable $w \in W(X')$, all
$E_i', F_j', G_k'$ are smooth rational curves.  
Consider the divisor $D =  2E_0' + E_1' + E_2' + E_3' + E_4'$.
Obviously $D^2 = 0$.  If necessary, by replacing $D$ by $w(D)$ where $w \in W(X')$, we may assume that 
$D$ defines an elliptic fibration.  Then $D$ is a singular fiber of type $I_0^*$ and $E_i'$ $(0\leq i \leq 4)$ are 
components of singular fibers.
Thus we may assume that $E_i'$ $(0\leq i \leq 4)$ are smooth rational curves.  Next consider the divisor
$D' = 2E_0' + E_1' + E_2' + E_3' + E_5'$.  By replacing $D'$ by $w(D'), \ w \in W(X')$ with $w(E_i') = E_i'$,
$0 \leq i \leq 4$, $D'$ defines an elliptic fibration.  Thus we may assume that all $E_i'$ are smooth rational curves.
Since $\mid F_i' + G_i' \mid = \mid 2E_0' + E_1' + E_2' + E_3' + E_4' + E_5' - E_i' \mid$, all $F_i', G_i'$ are also
smooth rational curves.

Next we shall show that the incidence relation of $E_i', F_j', G_k'$ is the same as that of $E_i, F_j, G_k$.
Obviously $E_0'$ is pointwisely fixed by
$\sigma'$.  Recall that $\sigma'$ acts on $H^0(X',\Omega_{X'})$ non trivially.  By considering the action
of $\sigma'$ on the tangent space of $E_i' \cap E_0'$, $\sigma'$ acts on $E_i'$ non trivially.  Now consider the elliptic
fibration defined by the linear system $\mid 2E_0' + E_1' + E_2' + E_3' + E_4' + E_5' - E_i' \mid$ with sections $F_j'$, $G_j'$
$(j \not= i)$.  Since no elliptic curves have 
an automorphism of order 5, $\sigma'$ acts on the sections $F_j'$ and $G_j'$ non trivially.  Note that $\sigma'$ has exactly
two fixed points on each of $E_i' , F_j', G_k', \ (1\leq i, j, k \leq 5)$.  Hence $F_i'$ and $G_i'$ meets at one point with
multiplicity 2.  Now we can easily see that
not only the dual graph, but also the incident relation of $E_i', F_j', G_k'$ coincides with that of $E_i, F_j, G_k$.

Finally define the isometry $\iota'$ of order 2 of $S_Y$ by $\iota'(F_i') = G_i'$ $(1 \leq i \leq 5)$ and $\iota'(E_i') = E_i'$.
Then $\iota'$ can be extended to an isometry of $H^2(X', {\bbZ})$ acting on $T_{X'}$ as $-1_{T_{X'}}$ as $\iota$.
By definition of $\iota'$, it preserves $C(M)$, and hence preserves the K\"ahler cone
$C(X)$ (Lemma \ref{V2}).
By the Torelli theorem, there exists an automorphism $\tau'$ with $(\tau')^* = \iota'$.  
It follows from Nikulin \cite{N3}, Theorem 4.2.2 that the set of fixed poins of $\tau'$ is the disjoint union of $E_0$ and
a smooth curve of genus 6.  By taking the quotient of $X'$ by $\tau'$, we have the same configuration as in 
\ref{K3}.  Thus $X'$ can be deformed to $X$ smoothly and hence $\sigma^*$ is conjugate to $\rho$.

\end{proof}

\subsection{Lemma}\label{compact}
{\it
Let $e \in T$ with $e^{2} = 0$.  Let $K$ be the sublattice generated by $\rho^{i}(e)$ $(0 \leq i \leq 4)$.
Then $K$ contains a vector with positive norm.
}

\begin{proof}
First note that $e$, $\rho(e)$, $\rho^{2}(e)$ are linearly independent isotoropic vectors.
Since the signature of $T$ is $(2,8)$, we may assume that $ \langle e, \pm \rho(e) \rangle > 0.$   
Then $e \pm \rho(e)$ is a desired one.
\end{proof}

\subsection{Lemma}\label{root}
{\it
Let $r \in T$  with $r^{2} = -2$.
Let $R$ be the lattice generated by  $\rho^{i}(r)$ $(0 \leq i \leq 4)$.
Assume that $R$ is negative definite.  Then
$R$ is isometric to the root lattice $A_{4}$.
}

\begin{proof}
Put $m_{i} = \langle r, \rho^{i}(r) \rangle$, $1 \leq i \leq 4$.
Then by assumption $\mid m_{i} \mid \leq 1$.
Also obviously $m_{1} = m_{4}, m_{2} = m_{3}$ and $\sum_{i=0}^{4} \rho^{i}(r) = 0$.
Then
$$-2 = r^{2} = \langle r, -\sum_{i=1}^{4} \rho^{i}(r)
\rangle = -2m_{1} -2m_{2}.$$
Hence $(m_{1}, m_{2} ) = (1,0)$ or $(0,1)$.
Therefore $\{ \rho^{i}(r) : 0 \leq i \leq 3 \}$ is a basis of the root lattice $A_{4}$.
\end{proof}

\subsection{Lemma}\label{root2}
{\it
Let $R \cong A_4$ be a sublattice of $T$.  Assume that $R$ is invariant under the action of $\rho$.
Then the orthogonal complement $R^{\perp}$ of $R$ in
$T$ is isomorphic to $U \oplus V \oplus A_4$.
}

\begin{proof}
Let $T'$ be the orthogonal complement of $R$ in $T$.  Then $T = R \oplus T'$ or
$T$ contains $R \oplus T'$ as a sublattice of index 5.  We shall show that the second case does not occur.
Assume that $[T : R \oplus T'] = 5$.  Then $A_{T'} = (T')^*/T' \cong (\bbF_5)^{\oplus 4}$ because
$\mid A_T\mid \cdot [T: R \oplus T']^2 = \mid A_R\mid \cdot \mid A_{T'}\mid$.
Let $\rho'$ be an isometry of $L$ so that $\rho' \mid T' = \rho \mid T'$ and $\rho' \mid (T')^{\perp} = 1$.  The existence
of such $\rho'$ follows from \cite{N1}, Corollary 1.5.2.
It follows from the surjectivity of the period map of $K3$ surfaces that there exists a $K3$ surface $Y$ whose
transcendental lattice isomorphic to $T'$ and whose period is an eigen-vector of $\rho'$ under a suitable
marking.  Since $\rho'$ acts trivially on the Picard lattice $(T')^{\perp}$ of $Y$, $\rho'$ is induced from an
automorphism $\sigma'$ of $Y$.  It follows from Vorontsov's theorem \cite{Vo} that the number of minimal
generator of $A_{T'}$ is at most ${\rm rank}(T') / \varphi (5) = 2$ where $\varphi$ is the Euler function.
This contradicts the fact $A_{T'} \cong (\bbF_5)^{\oplus 4}$.  Thus we have proved that
$T = R \oplus T'$.  Since $q_{T'} \cong q_{U \oplus V \oplus A_4}$, the assertion now follows from Nikulin
\cite{N1}, Theorem 1.14.2.
\end{proof}

\subsection{Discriminant locus}\label{root3}
Let $r \in T$ with $r^2 = -2$.  Let $R$ be the sublattice generated by  $\rho^{i}(r)$ $(0 \leq i \leq 4)$.
Assume that $R$ is negative definite.  Then $R \cong A_4$ and the orthogonal complement of $R$ in $T$ 
is isomorphic to $T' = U \oplus V \oplus A_4$ (Lemmas \ref{root}, \ref{root2}).  Let $\rho'$ be an isometry of $L$
so that $\rho' \mid T' = \rho \mid T'$ and $\rho' \mid (T')^{\perp} = 1$.  Then there exists an $K3$ surface
$Y$ and an automorphism $\sigma'$ such that the transcendental lattice of $Y$ is isomorphic to $T'$, 
the period of $Y$ is an eigen-vector of $\rho'$ and $\sigma'$ acts trivially on the Picard lattice of $Y$.
By the same argument as in the proof of Lemma \ref{isometry5}, we can see that $Y$ is corresponding to the case that
$f_5 = 0$ has a multiple root in \ref{K3}.

\section{A complex ball uniformization}\label{uniformization}

\subsection{Hermitian form}\label{hermite}
Let $\zeta = e^{4\pi \sqrt{-1}/5}$.
We consider $T$ as a free ${\bbZ}[\zeta]$-module $\Lambda$ by
$$(a+b\zeta)x = ax + b\rho(x).$$
Let 
$$h(x,y) = {2 \over 5 + \sqrt{5}} \sum_{i=0}^{4} \  \zeta^i \langle x, \rho^i(y) \rangle.$$
Then $h(x,y)$ is a hermitian form on ${\bbZ}[\zeta]$-module $\Lambda$.
With respect to a ${\bbZ}[\zeta]$-basis $e_1$ of $A_4$,
the hermitian matrix of $h\mid A_4$ is given by $-1$.
And with respect to a ${\bbZ}[\zeta]$-basis $e$ of $U\oplus V$,
the hermitian matrix of $h\mid U\oplus V$ is given by $(\sqrt{5} -1)/2$.
Thus $h$ is given by

\begin{equation}\label{hermitian}
\begin{pmatrix}{\sqrt{5} -1 \over 2}&0&0
\\0&-1&0
\\0&0&-1
\end{pmatrix}.
\end{equation}

\noindent
Let 
$$\varphi : \Lambda \to T^*$$
be a linear map defined by $$\varphi(x) = \sum_{i=0}^{3} \  (i+1)\rho^i(x)/5.$$
Note that $\varphi((1-\zeta)x) = \varphi(x - \rho(x)) = -\rho^4(x) \in T$. Hence
$\varphi$ induces an isomorphism
\begin{equation}\label{cong}
\Lambda/(1-\zeta)\Lambda \simeq A_T = T^*/T.
\end{equation}

\subsection{Reflections}\label{reflection}
Let $a \in \Lambda$ with $h(a,a) = -1$.  Then the map
$$R_a^{\pm} : v \to v  - (-1\pm \zeta) h(v, a) a$$
is an automorphism of $\Lambda$.  This automorphism $R_a^+$ has order 5 and $R_a^-$ has order 10
both of which fix the orthogonal complement of $a$.  They are called reflections.
Consider a decomposition 
$$T = U \oplus V \oplus A_4 \oplus A_4.$$
If $a = e_1$ of the last component $A_4$ as in \eqref{A4}, we can easily see that
$$R_a^{\pm} = \pm s_{e_1} \circ s_{e_2} \circ s_{e_3} \circ s_{e_4}$$
where $s_{e_i}$ is a reflection in $O(T)$ associated with $(-2)$-vector $e_i$ defined by
$$s_{e_i} : x \to x + \langle x, e_i \rangle e_i.$$
In other words, 
$$R_{e_1}^{\pm} = 1_U \oplus 1_V \oplus 1_{A_4} \oplus (\pm \rho_4).$$
Since $s_{e_i}$ acts trivially on $A_T$, $R_a^+$ acts trivially on $A_T \simeq \Lambda/(1-\zeta)\Lambda$.
On the other hand, $R_a^-$ acts on $A_T$ as a reflection associated with 
$\alpha = (e_1 + 2e_2 + 3e_3 + 4e_4)/5 \in A_T.$  

\subsection{The period domain and arithmetic subgroups}\label{period}
We use the same notation as in \ref{isometry}.
Let
$$T \otimes {\bbC} = T_{\zeta} \oplus T_{\zeta^{2}}
\oplus T_{\zeta^{3}} \oplus T_{\zeta^{4}}$$
be the decomposition of $\rho$-eigenspaces where $\zeta$ is a primitive
5-th root of unity (see Nikulin \cite{N2}, Theorem 3.1).  An easy calculation shows that 
$$\xi = e_1+(\zeta^4+1)e_2+(-\zeta-\zeta^2)e_3-\zeta e_4$$
is an eigenvector of $\rho_4$ with the eigenvalue $\zeta$ and
$$\langle \xi, \bar{\xi}\rangle = -5.$$
On the other hand, 
$$\mu = e-(\zeta^4+1)f+(-\zeta-\zeta^2)(e+f+y)-\zeta(-e+f-x)$$
is an eigenvector of $\rho_0$ with the eigenvalue $\zeta$ and 
$$\langle\mu, \bar{\mu}\rangle = 5(\zeta^2+\zeta^3).$$
Thus if $\zeta = e^{\pm 4\pi \sqrt{-1}/5}$, the hermitian form
$\langle \omega, {\bar \omega} \rangle /5$ on $T_{\zeta}$
is of signature $(1,2)$ and is given by
\begin{equation}\label{}
\begin{pmatrix}{\sqrt{5} -1 \over 2}&0&0
\\0&-1&0
\\0&0&-1
\end{pmatrix}.
\end{equation}
For other $\zeta$, the hermitian form is negative definite.
Now we take $\zeta = e^{4\pi \sqrt{-1}/5}$ and define

\begin{equation}\label{domain}
{\calB} = \{ z \in {\bbP}(T_{\zeta}) :
\langle z, {\bar z} \rangle > 0 \}.
\end{equation}
Then ${\calB}$ is a 2-dimensional complex ball.
For a $(-2)$-vector $r$ in $T$, we define
$${\calH}_r = r^{\perp} \bigcap {\calB}, \quad {\calH} = \bigcup_r {\calH}_r$$
where $r$ runs over $(-2)$-vectors in $T$.  
Let
\begin{equation}\label{arithmetic}
\Gamma = \{ \phi \in O(T) : \phi \circ \rho = \rho \circ \phi \},
\quad
\Gamma' = \{ \phi \in \Gamma : \phi \mid A_T = 1 \}.
\end{equation}

\subsection{Remark}\label{M-Y} 
The hermitian form $h$ in \eqref{hermitian} coincides with the one of Shimura \cite{S}, Yamazaki and Yoshida 
\cite{YY}.  This and the isomorphism \eqref{cong} imply that our groups $ \Gamma$, $\Gamma'$ coincide
with $ \Gamma$, $\Gamma(1-\mu)$ in Yamazaki and Yoshida \cite{YY}.

\subsection{Proposition}\label{prop6.2} 
{\it
$(1)$  $\Gamma$ is generated by reflections $R_a^-$ with $h(a,a) = -1$ and $\Gamma'$ is generated by
$R_a^+$ with $h(a,a) = -1$.  The quotient $\Gamma/ \Gamma'$ is isomorphic to
$\O(3,\bbF_5) \simeq {\bbZ}/2{\bbZ} \times S_5$.

$(2)$  ${\calH}/\Gamma'$ consists of $10$ smooth rational curves forming $10$ lines on the quintic
del Pezzo surface.
}

\begin{proof}
The assertions follow from the above Remark \ref{M-Y} and 
Propositions 4.2, 4.3, 4.4 in Yamazaki and Yoshida \cite{YY}.
\end{proof}

\subsection{Discriminant quadratic forms and discriminant locus}\label{}
Let
$$q_T : A_T \to \bbQ/2\bbZ$$
be the discriminant quadratic form of $T$.
The discriminant group $A_T$ consists of the following 125 vectors:
\smallskip

Type $(00): \alpha = 0, \ \# \alpha = 1$;

Type $(0): \alpha\not=0, \ q_T(\alpha) = 0, \ \# \alpha = 24$;

Type $(2/5): q_T(\alpha) = 2/5, \ \#\alpha = 30$;

Type $(-2/5):  q_T(\alpha) = -2/5, \ \#\alpha = 30$;

Type $(4/5):   q_T(\alpha) = 4/5, \ \#\alpha = 20$;

Type $(-4/5):  q_T(\alpha) = -4/5, \ \#\alpha = 20$.

\smallskip
\noindent
Let $A_4$ be a component of $T$ with a basis $e_1, e_2, e_3, e_4$ as in \ref{isometry}.
Then $(e_1 + 2e_2+3e_3+4e_4)/5 = (e_1 2\rho(e_1) + 3\rho^2(e_1) + 4\rho^3(e_1))/4$ mod $T$ is
a vector in $A_T$ with norm $-4/5$.  It follows from Proposition \ref{prop6.2} that $\Gamma/\Gamma'$ acts transitively
on the set of $(-4/5)$-vectors in $A_T$.  Hence
for each $\alpha \in A_T$ with $q_T(\alpha) = -4/5$ there exists a vector $r \in T$ with $r^2 =-2$ satisfying 
$\alpha = (r + 2\rho(r) + 3\rho^2(r) + 4\rho^3(r))/5$.
Moreover $\pm \alpha$ defines $$\calH_{\alpha} = \bigcup_r H_r$$
where $r$ moves over the set 
$$\{ r \in T : r^2 = -2, \alpha = (r + 2\rho(r) + 3\rho^2(r) + 4\rho^3(r))/5\ {\rm mod} \ T\}.$$
Thus the set 
$$\{ \alpha \in A_T : q_T(\alpha) = -4/5 \}/\pm 1$$ 
bijectively corresponds to the set of components of
$\calH/\Gamma'$.
Let
$$\tilde{\Gamma} = \{ \tilde{\phi} \in O(L) : \tilde{\phi} \circ \rho = \rho \circ \tilde{\phi} \}.$$

\subsection{Lemma}\label{surj} 
{\it
The restriction map $\tilde{\Gamma} \to {\Gamma}$ is surjective.
}

\begin{proof}
We use the same notation as in \ref{K3}.  The symmetry group $S_5$ of degree 5 naturally acts on the
set $\{ E_1,..., E_5 \}$ as permutations.  This action can be extended to the one on $S$.
Together with the action of $\iota$, the natural map 
$$O(S) \to O(q_S) \cong \{\pm 1\} \times S_5$$
is surjective.  Let $g \in \Gamma$.  Then the above implies that there exists an isometry $g'$ in $O(S)$
whose action on $A_S \cong A_T$ coincides with the one of $g$ on $A_T$.  Then it follows from Nikulin
\cite{N1}, Proposition 1.6.1 that the isometry $(g', g)$ of $S\oplus T$ can be extended to an isometry in
$\tilde{\Gamma}$ which is the desired one.
\end{proof}

\subsection{Period map}\label{}

We shall define an $S_5$-equivariant map
$$p : P_1^5 \to {\calB} / \Gamma'$$
called the {\it period map}.
Denote by $(P_1^5)^0$ the locus of distinct five ordered points on $\bbP^1$.
Let $\{p_1,...,p_5\} \in (P_1^5)^0$. 
Let $X$ be the corresponding $K3$ surface with the automorphism $\sigma$ of order 5 as in \ref{K3}.
The order of $\{p_1,...,p_5\}$ defines an order of smooth rational curves
$$E_i, \quad ( 0 \leq i \leq 5) \quad F_j, G_j, \quad ( 1 \leq j \leq 5)$$
modulo the action of the covering involution $\iota$.  
It follows from Lemma  \ref{isometry5} that 
there exists an isometry
$$\alpha : L \to H^{2}(X, {\bbZ})$$
satisfying $\alpha  \circ  \rho =  \sigma^* \circ  \alpha$.
Now we define 
$$p(X, \alpha) = (\alpha\otimes \bbC )^{-1}(\omega_X).$$

\subsection{Lemma}\label{} 
$p(X, \alpha) \in {\calB} \setminus \calH$.

\begin{proof}
If not, there exists a vector $r \in T$ with $r^2 = -2$ which is represented by an effective divisor on $X$ as in
the proof of Lemma \ref{-2}.  Obviously 
$r + \sigma^* (r) + \cdot \cdot \cdot + (\sigma^*)^4 (r) = 0$.  On the other hand
$r + \sigma^* (r) + \cdot \cdot \cdot + (\sigma^*)^4 (r)$ is non-zero effective
because $\sigma$ is an automorphism.  Thus we have a contradiction. 
\end{proof}

Thue we have a holomorphic map
$$p : (P_1^5)^0 \to (\calB \setminus \calH)/ \Gamma'.$$
The group $S_5$ naturally acts on $P_1^5$ which induces an action on $S$ as permutations of $E_1,..., E_5$.
On the other hand, $S_5\cong \Gamma/\{\pm 1\} \cdot \Gamma'$ naturally acts on $\calB/ \Gamma'$.
Under the natural isomorphism $O(q_S) \cong O(q_T)\cong \{\pm 1\}\cdot S_5$, $p$ is equivariant under these actions of $S_5$.

It is known that the quotient ${\calB}/ \Gamma'$ is compact (see Shimura \cite{S}).
We remark that cusps of $\calB$ correspond to totally isotropic sublattices of
$T$ invariant under $\rho$.  Hence the compactness also follows from Lemma \ref{compact}.

\subsection{Main theorem}\label{Main}
{\it
The period map $p$ can be extended to an $S_5$-equivariant isomorphism}
$$\tilde{p} : P_1^5 \to {\calB} / \Gamma'.$$

\begin{proof}
Let $\calM$ be the space of all 5 stable points on $\bbP^1$ and $\calM_0$ the space
of all distinct 5 points on $\bbP^1$.  We can easily see that $\calM \setminus \calM_0$ is locally
contained in a divisor with normal crossing.  By construction, $p$ is locally liftable to $\calB$.
It now follows from a theorem of Borel
\cite{Borel} that $p$ can be extended to a holomorphic map from $\calM$ to 
$\calB/\Gamma'$ which induces a holomorphic map  
$\tilde{p}$ from $P_1^5$ to $\calB /\Gamma'$.
Next we shall show the injectivity of the period map
over $(\calB \setminus \calH)/ \Gamma'$.  Let $C, C'$ be two plane quintic curves as in \ref{quintic}.
Let $(X, \alpha)$ (resp. $(X', \alpha')$) be the associated marked $K3$ surfaces with automorphisms $\tau, \sigma$
(resp. $\tau', \sigma'$).  Assume that the periods of $(X, \alpha)$ and $(X', \alpha')$ coincide in
$\calB/\Gamma'$.  Then there exists an isometry
$$\varphi : H^2(X', {\bbZ}) \to H^2(X, {\bbZ})$$
preserving the periods and satisfying
$\varphi \circ (\tau')^* = \tau^* \circ \varphi$ and $\varphi \circ (\sigma')^* = 
\sigma^* \circ \varphi$ (Lemma \ref{surj}).
It follows from Lemma \ref{V} that $\varphi$ preserves the K\"ahler cones.
The Torelli theorem for $K3$ surfaces implies that
there exists an isomorphism $f: X \to X'$ with $f^* = \varphi$.
Then $f$ induces an isomorphism between the corresponding plane quintic curves $C$ and $C'$.
Thus we have proved the injectivity of the period map.

Since both $P_1^5$ and $\calB / \Gamma'$ are compact,
$\tilde{p}$ is surjective.  Recall that both $P_1^5 \setminus (P_1^5)^0$ and 
$\calH /\Gamma'$ consist of 10 smooth rational curves.
The surjectivity of $\tilde{p}$ implies that no components of $P_1^5 \setminus (P_1^5)^0$ contract to a point.
Now the Zariski main theorem implies that $\tilde{p}$ is isomorphic.  
By construction, $\tilde{p}$ is $S_5$-equivariant
over the Zariski open set $(P_1^5)^0$.  Hence $\tilde{p}$ is $S_5$-equivariant isomorphism between
$P_1^5$ and $\calB /\Gamma'$.
\end{proof}

\section{Shimura-Terada-Deligne-Mostow's reflection groups}\label{Shimura}
The plane quintic curve $C$ defined by \eqref{quintic} 
appeared in the papers of Shimura \cite{S}, Terada \cite{Te}, 
Deligne-Mostow \cite{DM}, and the moduli space of these curves has a complex ball uniformization.  
As we remarked, the hermitian form \eqref{hermitian} coincides with
those of Shimura \cite{S}, Terada \cite{Te}, Deligne-Mostow \cite{DM}  
(see Remark \ref{M-Y}).   This implies

\subsection{Theorem}\label{Mostow}
{\it
The arithmetic subgroup $\Gamma$ is  the one appeared in
Deligne-Mostow's list \ \cite{DM}$:$}
$$\begin{pmatrix}{2 \over 5}, &{2 \over 5}, &{2 \over 5}, &{2 \over 5}, &{2 \over 5}\end{pmatrix}.$$

A geometric meaning of this theorem is as follows.
Recall that $X$ has an {\it isotrivial} pencil of curves of genus two whose
general member is the smooth curve $D$ of genus two with an automorphism
of order 5 given by the equation \eqref{bolza}.
The $X$ is given by
$$s^2 = x_0(x_0^5 - f_5(x_1,x_2))$$
where $x_1/x_2$ is the parameter of this pencil.  On the other hand, we consider $C$ as
a base change $C \to \bbP^1$ given by
$$(v, x_1, x_2) \to (x_1, x_2).$$
Then over $C$, $v^5 = f_5(x_1,x_2)$ and hence the pencil is given by
$$s^2 = x_0(x_0^5 + v^5)$$
which is nothing but the equation of the curve $D$.  Thus
the $K3$ surface $X$ is birational to the quotient of $C \times D$ by an diagonal action of ${\bbZ}/5\bbZ$.
This correspondence gives a relation between the Hodge structures of $C$ and $X$.

\subsection{Problem}\label{}
Let $\mu_i$ be a positive rational number $(0 \leq i \leq d+1$ or $i = \infty )$ satisfying
$\sum_i \mu_i = 2$. 
Set 
$$F_{gh}(x_2,...,x_{d+1}) = \int_g^h u^{-\mu_0} (u-1)^{-\mu_1} \prod^{d+1}_{i=2} (u-x_i)^{-\mu_i} du$$
where $g, h \in \{ \infty, 0, 1, x_2,..., x_{d+1} \}$.  Then $F_{gh}$ is a multivalued function on
$$M = \{ (x_i) \in (\bbP^1)^{d+3} \ | \ x_i \not= \infty , 0, 1, \ x_i \not= x_j \ (i \not= j)\}.$$
These functions generate a $(d+1)$-dimensional vector space which is invariant under
monodromy.  Let $\Gamma_{(\mu_i)}$ be the image of $\pi_1(M)$ in $\PGL(d+1, \bbC)$ under
the monodromy action.  In Deligne-Mostow \cite{DM} and Mostow \cite{Mo}, they gave 
a sufficient condition for which $\Gamma_{(\mu_i)}$ is a lattice in the projective unitary group 
$PU(d,1)$, that is, $\Gamma_{(\mu_i)}$ is discrete and of finite covolume, and gave a list of such $(\mu_i)$
(see \cite{Th} for the correction of their list).

Denote $\mu_i = \mu'_i/D$ where $D$ is the common denominator.  As remarked in Theorem \ref{Mostow},
in the case $D=5$, $\Gamma_{(\mu_i)}$  is related to $K3$ surfaces.  In case of $D=3,4$ or $6$, 
$\Gamma_{(\mu_i)}$ is also related to
$K3$ surfaces (see \cite{K2}, \cite{K3}, \cite{DGK}).  In these cases, the corresponding $K3$ surfaces have
an isotrivial elliptic fibration whose general fiber is an elliptic curve with an automorphism of order 4 or 6.

For the remaining arithmetic subgroups $\Gamma_{(\mu_i)}$ with $D> 6$ in the Deligne-Mostow's list, 
are they  related to $K3$ surfaces ?


\end{document}